\documentclass[11pt]{amsart}
\usepackage{amsmath, amsthm, amsfonts, amssymb} 
\usepackage{url}
\usepackage{mathrsfs}  %% if you have trouble compiling this, comment out this and the next line
\renewcommand{\mathcal}{\mathscr}

\newtheorem{theo}{Theorem}[section]
\newtheorem*{TheoremA}{Theorem A}
\newtheorem*{TheoremB}{Theorem B}
\newtheorem*{CorollaryA}{Corollary A}
\newtheorem*{CorollaryB}{Corollary B}

\newtheorem{lem}[theo]{Lemma}
\newtheorem{prop}[theo]{Proposition}
\newtheorem{claim}[theo]{Claim}
\theoremstyle{definition}
\newtheorem{rem}[theo]{Remark}
\newtheorem{exam}[theo]{Example}

\newtheorem{nota}[theo]{Notation}
\newtheorem{df}[theo]{Definition}

\newcommand{\R}{\mathbf{R}}
\newcommand{\C}{\mathbf{C}}
\newcommand{\Z}{\mathbf{Z}}
\newcommand{\F}{\mathbf{F}}

\newcommand{\N}{\mathbf{N}}

\newcommand{\Id}{\operatorname{Id}}
\newcommand{\Ad}{\operatorname{Ad}}
\newcommand{\id}{\operatorname{id}}

\newcommand{\SL}{\operatorname{SL}}

\newcommand{\Haar}{\operatorname{Haar}}

\newcommand{\supp}{\operatorname{supp}}
\newcommand{\op}{\operatorname{op}}

\newcommand{\NC}{\operatorname{NC}}

\newcommand{\cb}{\operatorname{cb}}

\newcommand{\Lim}{\operatorname{Lim}}

\newcommand{\bin}{\operatorname{bin}}

\begin{document}

\title[strongly solid II$_1$ factors with an exotic MASA]{strongly solid II$_1$ factors \\ with an exotic MASA}

\begin{abstract}
Using an extension of techniques of Ozawa and Popa, we give an example of a non-amenable strongly solid $\rm{II}_1$ factor $M$ containing an ``exotic'' maximal abelian subalgebra $A$: as an $A$,$A$-bimodule, $L^2(M)$ is neither coarse nor discrete.  Thus we show that there exist $\rm{II}_1$ factors with such property but without Cartan subalgebras.  It also follows from Voiculescu's free entropy results that $M$ is not an interpolated free group factor, yet it is strongly solid and has both the Haagerup property and the complete metric approximation property.
\end{abstract}

\author[C. Houdayer]{Cyril Houdayer}
\address{CNRS ENS Lyon \\
UMPA UMR 5669 \\
69364 Lyon cedex 7 \\
France}
\email{cyril.houdayer@umpa.ens-lyon.fr}

\author[D. Shlyakhtenko]{Dimitri Shlyakhtenko*}
\address{UCLA\\	
Department of Mathematics\\
520 Portola Plaza\\
LA\\ 
CA 90095}
\email{shlyakht@math.ucla.edu}
\thanks{* Research partially supported by NSF grant DMS-0555680}


\maketitle

\section{Introduction}

In their breakthrough paper \cite{ozawapopa}, Ozawa and Popa showed that the free group factors $L(\F_n)$ are {\em strongly solid}, i.e. the normalizer $\mathcal{N}_{L(\F_n)}(P)=\{u\in \mathcal{U}(L(\F_n)): uPu^*=P\}$ of any diffuse amenable subalgebra $P\subset L(\F_n)$ generates an amenable von Neumann algebra, thus AFD by Connes' result \cite{connes76}. This strengthened two well-known indecomposability results for free group factors: Voiculescu's celebrated result in \cite{voiculescu96}, showing that $L(\F_n)$ has no Cartan subalgebra, which in fact exhibited the first examples of factors with no Cartan decomposition; and Ozawa's result in \cite{ozawa2003}, showing that the commutant in $L(\F_n)$ of any diffuse subalgebra must be amenable ($L(\F_n)$ are {\it solid}). Furthermore in \cite{ozawapopaII}, Ozawa and Popa showed that for any lattice $\Gamma$ in $\SL(2, \R)$ or $\SL(2, \C)$, the group von Neumann algebra $L(\Gamma)$ is strongly solid as well.

In this paper, we use a combination of Popa's deformation and intertwining techniques \cite{{popasup}, {popamal1}, {popa2001}} and the techniques of Ozawa and Popa \cite{ozawapopa, ozawapopaII} to give another example of a strongly solid  $\rm{II_1}$ factor not isomorphic to an amplification of a free group factor, i.e. to an {\em interpolated} free group factor \cite{{dykema94}, {radulescu1994}} (the first example of this kind was constructed by the first-named author in \cite{houdayer7}, answering an open question of Popa \cite{popa07}).

Our example is rather canonical: it is the crossed product of a free group factor $L(\F_\infty)$ by $\mathbf{Z}$, acting by a free Bogoljubov transformation obtained via Voiculescu's free Gaussian functor (cf. \cite{voiculescu92}). Roughly speaking, recall \cite{voiculescu92} that to any separable real Hilbert space $H_\R$, one can associate a finite von Neumann algebra $\Gamma(H_\R)''$ which is precisely isomorphic to the free group factor $L(\F_{\dim H_\R})$. To any orthogonal representation $\pi : \Z \to \mathcal{O}(H_\R)$ of $\Z$ on $H_\R$ corresponds a trace-preserving action $\sigma^\pi : \Z \curvearrowright \Gamma(H_\R)''$, called the {\em Bogoljubov action} associated with the orthogonal representation $\pi$. Alternatively, our algebra can be viewed as a free Krieger algebra in the terminology of \cite{shlya99}, constructed from an abelian subalgebra and a certain completely positive map (related to the spectral measure of the $\mathbf{Z}$-action). It is in this way rather similar to a core of a free Araki-Woods factor  \cite{shlya98,shlya97}. Along these lines, our main results are the following.

\begin{TheoremA}
Let $\pi : \Z \to \mathcal{O}(H_\R)$ be an orthogonal representation on the real Hilbert space $H_\R$ such that the spectral measure of $\pi$ has no atoms. Denote by $M = L(\F_\infty) \rtimes_{\sigma^\pi} \Z$ the crossed product under the Bogoljubov action. Then for any maximal abelian subalgebra $A \subset M$, the normalizer $\mathcal{N}_M(A)$ generates an amenable von Neumann algebra. 
\end{TheoremA}

In particular, the ${\rm II_1}$ factor $M = L(\F_\infty) \rtimes_{\sigma^\pi} \Z$ has no Cartan subalgebras. Under additional assumptions on the orthogonal representation $\pi$, we can obtain a stronger result.

\begin{TheoremB}
Let $\pi : \Z \to \mathcal{O}(H_\R)$ be a {\em mixing} orthogonal representation on the real Hilbert space $H_\R$. Then $M = L(\F_\infty) \rtimes_{\sigma^\pi} \Z$ is a non-amenable strongly solid ${\rm II_1}$ factor, i.e. for any $P \subset M$ diffuse amenable subalgebra, $\mathcal{N}_M(P)''$ is an amenable von Neumann algebra.
\end{TheoremB}

Note that in both cases, $M$ has the Haagerup property and the complete metric approximation property, i.e. $\Lambda_{\cb}(M) = 1$.

The proof of Theorems A and B, following a {\textquotedblleft deformation/rigidity\textquotedblright} strategy, is a  combination of the ideas and techniques in \cite{{houdayer7}, {ozawapopa}, {ozawapopaII}, {popamal1}}. We will use the {\textquotedblleft free malleable deformation\textquotedblright} by automorphisms $(\alpha_t, \beta)$ defined on $\Gamma(H_\R)'' \ast \Gamma(H_\R)'' = \Gamma(H_\R \oplus H_\R)''$. This deformation naturally arises as the {\textquotedblleft second quantization\textquotedblright} of the rotations/reflection defined on $H_\R \oplus H_\R$ that commute with the $\Z$-representation $\pi \oplus \pi$.

The proof of Theorem B then consists in two parts. Let $\pi : \Z \to \mathcal{O}(H_\R)$ be a mixing orthogonal representation and denote by $M = L(\F_\infty) \rtimes_{\sigma^\pi} \Z$ the corresponding crossed product ${\rm II_1}$ factor. First, we show that given any amenable subalgebra $P \subset M$ such that $P$ does not embed into $L(\Z)$ inside $M$, the normalizer $\mathcal{N}_M(P)$ generates an amenable von Neumann algebra (see Theorem \ref{step}). For this, we will exploit the facts that the deformation $(\alpha_t)$ does not converge uniformly on the unit ball $(P)_1$ and that $P \subset M$ is {\it weakly compact}, and use the technology from \cite{{ozawapopa}, {ozawapopaII}}. So if $P \subset M$ is diffuse, amenable such that $\mathcal{N}_M(P)''$ is not amenable, $P$ must embed into $L(\Z)$ inside $M$. Exploiting Popa's intertwining techniques and the fact that the $\Z$-action $\sigma^\pi$ is  mixing, we prove that $\mathcal{N}_M(P)''$ is {\textquotedblleft captured\textquotedblright} in $L(\Z)$ and finally get a contradiction.

In proving that free group factors $L(\F_n)$ have no Cartan subalgebras \cite{voiculescu96}, Voiculescu
proved that they actually have a formally stronger property: for any MASA (maximal abelian subalgebra) $A\subset N=L(\F_n)$, $L^2(N)$ (when viewed as an $A$,$A$-bimodule) contains a sub-bimodule of $L^2(A) \otimes L^2(A)$.
In more classical language, for every MASA $L^\infty[0,1]\cong A\subset N$, every vector $\xi \in L^2(N)$ gives rise to a measure $\psi=\psi_\xi$ on $[0,1]^2$ determined by $$\int f(x) g(y) d\psi (x,y) = \langle f Jg^*J \xi , \xi\rangle,\qquad f,g\in A.$$ Voiculescu proved that, for any such $A\subset N\cong L(\F_n)$, there exists a nonzero vector $\xi$ for which $\psi$ is Lebesgue absolutely continuous. Any $N$ with this property cannot of course have Cartan subalgebras, since if $A$ is a Cartan subalgebra, the measure $\psi$ will have to be ``$r$-discrete'' (i.e., $\psi (B) = \int \nu_t(B) dt$ for some family of discrete measures $\nu_t$).

This raised the obvious question: if $N$ has no Cartan subalgebras, must it be that for any diffuse MASA $A\subset N$, the $A$,$A$-bimodule $L^2(N)$ contains a sub-bimodule of $L^2(A) \otimes L^2(A)$?
We answer this question in the negative.  Our examples $M = L(\F_\infty) \rtimes \Z$, while strongly solid (or having no Cartan subalgebras), have an ``exotic'' MASA $A = L(\Z)$, so that $L^2(M)$, when viewed as an $A$,$A$-bimodule, contains neither coarse nor $r$-discrete sub-bimodules.  In other words, for all $\xi\neq 0$, $\psi_\xi$ is neither $r$-discrete nor Lebesgue absolutely continuous. In particular, combined with Voiculescu's results, this property shows that our examples $M$ are not interpolated free group factors. Thus we prove:

\begin{CorollaryA}
Let $\pi : \Z \to \mathcal{O}(H_\R)$ be an orthogonal representation on the real Hilbert space $H_\R$ such that the spectral measure of $\bigoplus_{n \geq 1} \pi^{\otimes n}$ is singular w.r.t. the Lebesgue measure and has no atoms. Then the non-amenable ${\rm II_1}$ factor $M = L(\F_\infty) \rtimes_{\sigma^\pi} \Z$ has no Cartan subalgebra and is not isomorphic to any interpolated free group factor $L(\F_t)$, $1 < t \leq +\infty$.
\end{CorollaryA}

Assuming that the representation $\pi$ is mixing, we can obtain (see Theorem \ref{singular-ssolid}) new examples of strongly solid ${\rm II_1}$ factors not isomorphic to interpolated free group factors (see \cite{{houdayer7}, {popa07}}).

\begin{CorollaryB}
Let $\pi : \Z \to \mathcal{O}(H_\R)$ be a mixing orthogonal representation on the real Hilbert space $H_\R$ such that the spectral measure of $\bigoplus_{n \geq 1} \pi^{\otimes n}$ is singular w.r.t. the Lebesgue measure. Then the non-amenable ${\rm II_1}$ factor $M = L(\F_\infty) \rtimes_{\sigma^\pi} \Z$ is strongly solid and is not isomorphic to any interpolated free group factor $L(\F_t)$, $1 < t \leq +\infty$.
\end{CorollaryB}

In Section \ref{examples}, we will present examples of orthogonal representations $\pi : \Z \to \mathcal{O}(H_\R)$ which satisfy the assumptions of Corollaries A and B. After recalling the necessary  background in Section \ref{preliminaries}, Theorems A and B are proven in Section \ref{keyresult}.

{\bf Acknowledgements.} Part of this work was done while the first-named author was at University of California, Los Angeles. The first-named author is very grateful to the warm hospitality and the stimulating atmosphere at UCLA. He finally thanks Stefaan Vaes for fruitful discussions regarding this work during his visit at University of Leuven.

\section{Preliminaries}\label{preliminaries}

\subsection{Popa's intertwining techniques}

We first recall some notation. Let $P \subset M$ be an inclusion of finite von Neumann algebras. The {\it normalizer of} $P$ {\it inside} $M$ is defined as
\begin{equation*}
\mathcal{N}_M(P) := \left\{ u \in \mathcal{U}(M) : \Ad(u) P = P \right\},
\end{equation*}
where $\Ad(u) = u \cdot u^*$. The inclusion $P \subset M$ is said to be {\it regular} if $\mathcal{N}_M(P)'' = M$. The {\it quasi-normalizer of} $P$ {\it inside} $M$ is defined as
\begin{equation*}
\mathcal{QN}_M(P) := \left\{ a \in M : \exists b_1, \dots, b_n \in M, aP \subset \sum_i Pb_i, Pa \subset \sum_i b_iP \right\}.
\end{equation*}
The inclusion $P \subset M$ is said to be {\it quasi-regular} if $\mathcal{QN}_M(P)'' = M$. Moreover,
\begin{equation*}
P' \cap M \subset \mathcal{N}_M(P)'' \subset \mathcal{QN}_M(P)''.
\end{equation*}
Let $A, B$ be finite von Neumann algebras. An $A, B$-{\it bimodule} $H$ is a complex (separable) Hilbert space $H$ together with two {\it commuting} normal $\ast$-representations $\pi_A : A \to \mathbf{B}(H)$, $\pi_B : B^{\op} \to \mathbf{B}(H)$. We shall intuitively write $a \xi b = \pi_A(x)\pi_B(y^{\op}) \xi$, $\forall x \in A, \forall y \in B, \forall \xi \in H$. We say that $H_B$ is {\it finitely generated} as a right $B$-module if $H_B$ is of the form $p L^2(B)^{\oplus n}$ for some projection $p \in \mathbf{M}_n(\C) \otimes B$.

In \cite{{popamal1}, {popa2001}}, Popa introduced a powerful tool to prove the unitary conjugacy of two von Neumann subalgebras of a tracial von Neumann algebra $(M, \tau)$. We will make intensively use of this technique. If $A, B \subset (M, \tau)$ are (possibly non-unital) von Neumann subalgebras, denote by $1_A$ (resp. $1_B$) the unit of $A$ (resp. $B$).

\begin{theo}[Popa, \cite{{popamal1}, {popa2001}}]\label{intertwining1}
Let $(M, \tau)$ be a finite von Neumann algebra. Let $A, B \subset M$ be possibly non-unital von Neumann subalgebras. The following are equivalent:
\begin{enumerate}
\item There exist $n \geq 1$, a possibly non-unital $\ast$-homomorphism $\psi : A \to \mathbf{M}_n(\C) \otimes B$ and a non-zero partial isometry $v \in \mathbf{M}_{1, n}(\C) \otimes 1_AM1_B$ such that $x v = v \psi(x)$, for any $x \in A$.

\item The bimodule $\vphantom{}_AL^2(1_AM1_B)_B$ contains a non-zero sub-bimodule $\vphantom{}_AH_B$ which is finitely generated as a right $B$-module. 

\item There is no sequence of unitaries $(u_k)$ in $A$ such that 
\begin{equation*}
\lim_{k \to \infty} \|E_B(a^* u_k b)\|_2 = 0, \forall a, b \in 1_A M 1_B.
\end{equation*}
\end{enumerate}
\end{theo}
If one of the previous equivalent conditions is satisfied, we shall say that $A$ {\it embeds into} $B$ {\it inside} $M$ and denote $A \preceq_M B$. For simplicity, we shall write $M^n := \mathbf{M}_n(\C) \otimes M$.

\subsection{The complete metric approximation property}

\begin{df}[Haagerup, \cite{haagerup79}]
A finite von Neumann algebra $(M, \tau)$ is said to have the {\it complete metric approximation property} (c.m.a.p.) if there exists a net $\Phi_n : M \to M$ of ($\tau$-preserving) normal finite rank completely bounded maps such that
\begin{enumerate}
\item $\lim_n \|\Phi_n(x) - x\|_2 = 0$, $\forall x \in M$;
\item $\lim_n \|\Phi_n\|_{\cb} = 1$.
\end{enumerate}
\end{df}
It follows from Theorem $4.9$ in \cite{anan95} that if $G$ is a countable amenable group and $Q$ is a finite von Neumann algebra with the c.m.a.p., then for any action $G \curvearrowright (Q, \tau)$, the crossed product $Q \rtimes G$ has the c.m.a.p. as well. The notation $\bar{\otimes}$ will be used for the {\it spatial} tensor product.

\begin{df}[Ozawa \& Popa, \cite{ozawapopa}]
Let $\Gamma$ be a discrete group, let $(P, \tau)$ be a finite von Neumann algebra and let $\sigma : \Gamma \curvearrowright P$ be a $\tau$-preserving action. The action is said to be {\it weakly compact} if there exists a net $(\eta_n)$ of unit vectors in $L^2(P \bar{\otimes} \bar{P})_+$ such that
\begin{enumerate}
\item $\lim_n \|\eta_n - (v \otimes \bar{v})\eta_n\|_2 = 0$, $\forall v \in \mathcal{U}(P)$; 
\item $\lim_n \|\eta_n - (\sigma_g \otimes \bar{\sigma}_g)\eta_n\|_2 = 0$, $\forall g \in \Gamma$;
\item $\langle (a \otimes 1)\eta_n, \eta_n \rangle = \tau(a) = \langle \eta_n, (1 \otimes \bar{a}) \eta_n \rangle$, $\forall a \in M, \forall n$.
\end{enumerate}
These conditions force $P$ to be amenable. A von Neumann algebra $P \subset M$ is said to be {\it weakly compact inside} $M$ if the action by conjugation $\mathcal{N}_M(P) \curvearrowright P$ is weakly compact.
\end{df}

\begin{theo}[Ozawa \& Popa, \cite{ozawapopa}]\label{weakcompact}
Let $M$ be a finite von Neumann algebra with the complete metric approximation property. Let $P \subset M$ be an amenable von Neumann subalgebra. Then $P$ is weakly compact inside $M$.
\end{theo}

\subsection{Voiculescu's free Gaussian functor \cite{DVV:free,voiculescu92}}

Let $H_\R$ be a real separable Hilbert space. Let $H = H_\R \otimes_\R \C$ be the corresponding complexified Hilbert space. The \emph{full Fock space} of $H$ is defined by
\begin{equation*}
\mathcal{F}(H) =\C\Omega \oplus \bigoplus_{n = 1}^{\infty} H^{\otimes n}.
\end{equation*}
The unit vector $\Omega$ is called the \emph{vacuum vector}. For any $\xi \in H$, we have the \emph{left creation operator}
\begin{equation*}
\ell(\xi) : \mathcal{F}(H) \to \mathcal{F}(H) : \left\{ 
{\begin{array}{l} \ell(\xi)\Omega = \xi, \\ \ell(\xi)(\xi_1 \otimes \cdots \otimes \xi_n) = \xi \otimes \xi_1 \otimes \cdots \otimes \xi_n.
\end{array}} \right.
\end{equation*}
For any $\xi \in H$, we denote by $s(\xi)$ the real part of $\ell(\xi)$ given by
\begin{equation*}
s(\xi) = \frac{\ell(\xi) + \ell(\xi)^*}{2}.
\end{equation*}
The crucial result of Voiculescu \cite{voiculescu92} is that the distribution of the operator $s(\xi)$ w.r.t. the vacuum vector state $\langle \cdot \Omega, \Omega\rangle$ is the semicircular law supported on the interval $[-\|\xi\|, \|\xi\|]$, and for any subset $\Xi\subset H_\R$ of pairwise orthogonal vectors, the family $\{s(\xi):\xi\in \Xi\}$ is freely independent. Set
\begin{equation*}
\Gamma(H_\R)'' = \{s(\xi) : \xi \in H_\R\}''.
\end{equation*}
The vector state $\tau = \langle \cdot \Omega, \Omega\rangle$ is a faithful normal trace on $\Gamma(H_\R)''$, and 
\begin{equation*}
\Gamma(H_\R)'' \cong L(\F_{\dim H_\R}).
\end{equation*}
Since $\Gamma(H_\R)''$ is a free group factor, $\Gamma(H_\R)''$ has the Haagerup property and the c.m.a.p. \cite{haagerup79}.

\begin{rem}[\cite{speicher:noncrossing, voiculescu92}]\label{value} %This formula is needed for the mixing property, Prop 2.8.
Explicitely the value of $\tau$ on a word in $s(\xi_\iota)$ is given by
\begin{equation}\label{formula}
\tau(s(\xi_1) \cdots s(\xi_n)) = 2^{-n}\sum_
{(\{\beta_i, \gamma_i\}) 
 \in \NC(n), 
\beta_i < \gamma_i}
\prod_{k = 1}^{n/2}\langle \xi_{\beta_k}, \xi_{\gamma_k}\rangle.
\end{equation}
for $n$ even and is zero otherwise. Here $\NC(2p)$ stands for all the non-crossing pairings of the set $\{1, \dots, 2p\}$, i.e. pairings for which whenever $a < b < c < d$, and $a, c$ are in the same class, then $b, d$ are not in the same class. The total number of such pairings is given by the $p$-th Catalan number
\begin{equation*}
C_p = \frac{1}{p + 1}\begin{pmatrix}
2p \\
p
\end{pmatrix}.
\end{equation*}
\end{rem}

Let $G$ be a countable group together with an orthogonal representation $\pi : G \to \mathcal{O}(H_\R)$. We shall still denote by $\pi : G \to \mathcal{U}(H)$ the corresponding unitary representation on the complexified Hilbert space $H = H_\R \otimes_\R \C$. The {\it free Bogoljubov shift} $\sigma^\pi : G \curvearrowright (\Gamma(H_\R)'', \tau)$ associated with the representation $\pi$ is defined by
\begin{equation*}
\sigma_g^\pi = \Ad(\mathcal{F}(\pi_g)), \forall g \in G,
\end{equation*}
where $\mathcal{F}(\pi_g) = \bigoplus_{n \geq 0} \pi_g^{\otimes n} \in \mathcal{U}(\mathcal{F}(H))$.
\begin{nota}
For a countable group $G$ together with an orthogonal representation $\pi : G \to \mathcal{O}(H_\R)$, we shall denote by
\begin{equation*}
\Gamma(H_\R, G, \pi)'' = \Gamma(H_\R)'' \rtimes_{\sigma^\pi} G.
\end{equation*}
\end{nota}

\begin{exam}
If $(\pi, H) = (\lambda_G, \ell^2(G))$ is the left regular representation of $G$, it is easy to see that the action $\sigma^{\lambda_G} : G \curvearrowright \Gamma(\ell^2(G))''$ is the {\it free} Bernoulli shift and in that case $\Gamma(\ell^2(G), G, \lambda_G)'' \cong L(\Z) \ast L(G)$.
\end{exam}

For any $n \geq 0$, denote by $K_\pi^{(n)} = H^{\otimes n} \otimes \ell^2(G)$ with the $L(G),L(G)$-bimodule structure given by:
\begin{eqnarray*}
u_g \cdot (\xi_1 \otimes \cdots \otimes \xi_n \otimes \delta_h) & = & \pi_g \xi_1 \otimes \cdots \otimes \pi_g \xi_n \otimes \delta_{gh} \\
(\xi_1 \otimes \cdots \otimes \xi_n \otimes \delta_h) \cdot u_g & = & \xi_1 \otimes \cdots \otimes \xi_n \otimes \delta_{hg}.
\end{eqnarray*}
It is then straightforward to check that as $L(G), L(G)$-bimodules, we have the following isomorphism
\begin{equation*}
L^2(\Gamma(H_\R, G, \pi)'') \cong \bigoplus_{n \geq 0} K_\pi^{(n)}.
\end{equation*}

Recall that $\pi$ is said to be {\it mixing} if 
\begin{equation*}
\lim_{g \to \infty} \langle \pi_g \xi, \eta\rangle = 0, \forall \xi, \eta \in H.
\end{equation*}
The following proposition is an easy consequence of Remark \ref{value} and Kaplansky density theorem.
\begin{prop}
Let $G$ be a countable group together with an orthogonal representation $\pi : G \to \mathcal{O}(H_\R)$. The following are equivalent:
\begin{enumerate}
\item The representation $\pi : G \to \mathcal{O}(H_\R)$ is mixing.
\item The $\tau$-preserving action $\sigma^\pi : G \curvearrowright \Gamma(H_\R)''$ is mixing, i.e. 
\begin{equation*}
\lim_{g \to \infty} \tau(\sigma_g^\pi(x)y) = 0, \forall x, y \in \Gamma(H_\R)'' \ominus \C.
\end{equation*}
\end{enumerate}
\end{prop}

\section{Proof of Theorems A and B}\label{keyresult}

\subsection{The free malleable deformation on $\Gamma(H_\R, G, \pi)''$}

Let $G$ be a countable group together with an orthogonal representation $\pi : G \to \mathcal{O}(H_\R)$. Set
\begin{itemize}
\item $M = \Gamma(H_\R, G, \pi)''$.
\item $\widetilde{M} = \Gamma(H_\R \oplus H_\R, G, \pi \oplus \pi)''$.
\end{itemize}
Thus, we can regard $\widetilde{M}$ as the amalgamated free product
\begin{equation*}
\widetilde{M} = M \ast_{L(G)} M,
\end{equation*}
where we view $M \subset \widetilde{M}$ under the identification with the left copy. Consider the following orthogonal transformations on $H_\R \oplus H_\R$:
\begin{eqnarray*}
U_t & = & 
\begin{pmatrix}
\cos(\frac{\pi}{2} t) & -\sin(\frac{\pi}{2} t) \\
\sin(\frac{\pi}{2} t) & \cos(\frac{\pi}{2} t)
\end{pmatrix}, \forall t \in \R, \\
V & = & 
\begin{pmatrix}
1 & 0 \\
0 & -1
\end{pmatrix}.
\end{eqnarray*}
Define the associated deformation $(\alpha_t, \beta)$ on $\Gamma(H_\R \oplus H_\R)''$ by 
\begin{equation*}
\alpha_t = \Ad(\mathcal{F}(U_t)), \; \beta = \Ad(\mathcal{F}(V)). 
\end{equation*}
Since $U_t, V$ commute with $\pi \oplus \pi$, it follows that $\alpha_t, \beta$ commute with the diagonal action $\sigma^\pi \ast \sigma^\pi$. We can then extend the deformation $(\alpha_t, \beta)$ to $\widetilde{M}$ by ${\alpha_t}_{|L(G)} = \beta_{|L(G)} = \Id$. Moreover it is easy to check that the deformation $(\alpha_t, \beta)$ is {\it malleable} in the sense of Popa:

\begin{prop}
The deformation $(\alpha_t, \beta)$ satisfies:
\begin{enumerate}
\item $\lim_{t \to 0} \|x - \alpha_t(x)\|_2 = 0$, $\forall x \in \widetilde{M}$.
\item $\beta^2 = \Id$, $\alpha_t \beta = \beta \alpha_{-t}$, $\forall t \in \R$.
\item $\alpha_1(x \ast_{L(G)} 1) = 1 \ast_{L(G)} x$, $\forall x \in M$.
\end{enumerate}
\end{prop}

We recall at last that the s-malleable deformation $(\alpha_t, \beta)$ automatically features a certain {\it transversality} property.
\begin{prop}[Popa, \cite{popasup}]\label{transversality}
We keep the same notation as before. We have the following:
\begin{equation}\label{trans}
\|x - \alpha_{2t}(x)\|_2 \leq 2 \|\alpha_t(x) - (E_{M} \circ \alpha_t)(x)\|_2, \; \forall x \in M, \forall t > 0.
\end{equation}
\end{prop}

The following result of the first-named author about intertwining subalgebras inside the von Neumann algebras $\Gamma(H_\R, G, \pi)''$ (see Theorems $5.2$ in \cite{houdayer3} and $3.4$ in \cite{houdayer6}) will be a crucial tool in the next subsection.

\begin{theo}[\cite{{houdayer3}, {houdayer6}}]\label{intertwining}
Let $G$ be a countable group. Let $\pi : G \to \mathcal{O}(H_\R)$ be any orthogonal representation. Set  $M = \Gamma(H_\R, G, \pi)''$. Let $p \in M$ be a non-zero projection. Let  $P \subset pMp$ be a von Neumann subalgebra such that the deformation $(\alpha_t)$ converges uniformly on the unit ball $(P)_1$. Then $P \preceq_M L(G)$.
\end{theo}

\subsection{The key result}

Let $M, N, P$ be finite von Neumann algebras. For any $M, N$-bimodules $H, K$, denote by $\pi_H$ (resp. $\pi_K$) the associated $\ast$-representation of the binormal tensor product $M \otimes_{\bin} N^{\op}$ on $H$ (resp. on $K$). We refer to \cite{EL} for the definition of $\otimes_{\bin}$. We say that $H$ is {\em weakly contained} in $K$ and denote it by $H \prec K$ if the representation $\pi_H$ is weakly contained in the representation $\pi_K$, that is if $\ker(\pi_H) \supset \ker(\pi_K)$. Let $H, K$ be $M, N$-bimodules. The following are true:

\begin{enumerate}
\item Assume that $H \prec K$. Then, for any $N$-$P$ bimodule $L$, we have $H \otimes_N L \prec K \otimes_N L$, as $M, P$-bimodules. Exactly in the same way, for any $P, M$-bimodule $L$, we have $L \otimes_M H \prec L \otimes_M K$, as $P, N$-bimodules (see Lemma $1.7$ in \cite{anan95}).

\item A von Neumann algebra $B$ is amenable iff $L^2(B) \prec L^2(B) \otimes L^2(B)$, as $B$-$B$ bimodules.
\end{enumerate}

Let $B, M, N$ be von Neumann algebras such that $B$ is amenable. Let $H$ be any $M, B$-bimodule and let $K$ be any $B, N$-bimodule. Then, as $M, N$-bimodules, we have $H \otimes_B K \prec H \otimes K$ (straightforward consequence of $(1)$ and $(2)$).

\begin{lem}\label{weakcontainment}
Let $G$ be an amenable group together with an orthogonal representation $\pi : G \to \mathcal{O}(H_\R)$. Let $M = \Gamma(H_\R, G, \pi)''$. The $M, M$-bimodule $\mathcal{H} = L^2(\widetilde{M}) \ominus L^2(M)$ is weakly contained in the coarse bimodule $L^2(M) \otimes L^2(M)$. In particular, the left $M$-action on $\mathcal{H}$ extends to a u.c.p. map $\Psi : \mathbf{B}(L^2(M)) \to \mathbf{B}(\mathcal{H})$ whose range commutes with the right $M$-action.
\end{lem}

\begin{proof}
Set $B = L(G)$ which is amenable by assumption. By definition of the amalgamated free product $\widetilde{M} = M \ast_{L(G)} M$ (see \cite{voiculescu92}), we have as $M, M$-bimodules
\begin{equation*}
L^2(\widetilde{M}) \ominus L^2(M) \cong \bigoplus_{n \geq 1} \mathcal{H}_n,
\end{equation*} 
where 
\begin{equation*}
\mathcal{H}_n = L^2(M) \otimes_B \mathop{\overbrace{(L^2(M) \ominus L^2(B)) \otimes_B \cdots  \otimes_B (L^2(M) \ominus L^2(B))}}^{2n - 1} \otimes_B L^2(M).
\end{equation*}
Since $B = L(G)$ is amenable, the identity bimodule $L^2(B)$ is weakly contained in the coarse bimodule $L^2(B) \otimes L^2(B)$. From the standard properties of composition and weak containment of bimodules (see Lemma $1.7$ in \cite{anan95}), it follows that as $M, M$-bimodules
\begin{equation*}
\mathcal{H}_n \prec L^2(M) \otimes \mathop{\overbrace{(L^2(M) \ominus L^2(B)) \otimes \cdots \otimes (L^2(M) \ominus L^2(B))}}^{2n - 1} \otimes L^2(M).
\end{equation*}
Consequently, we obtain as $M, M$-bimodules
\begin{equation*}
\mathcal{H} = L^2(\widetilde{M}) \ominus L^2(M) \prec L^2(M) \otimes L^2(M).
\end{equation*}
Now the rest of the proof is the same as the one of Lemma $5.1$  in \cite{ozawapopaII}. The binormal representation $\mu$ of $M \odot M^{\op}$ on $\mathcal{H}$ is continuous w.r.t. the minimal tensor product. Hence $\mu$ extends to a u.c.p. map $\tilde{\mu}$ from $\mathbf{B}(L^2(M)) \bar{\otimes} M^{\op}$ to $\mathbf{B}(\mathcal{H})$. Define $\Psi(x) = \tilde{\mu}(x \otimes 1)$, $\forall x \in \mathbf{B}(L^2(M))$. Since $M^{\op}$ is in the multiplicative domain of $\tilde{\mu}$, it follows that the range of $\Psi$ commutes with the right $M$-action.
\end{proof}

The next theorem, which is the key result of this section in order to prove Theorems A and B, can be viewed as an analog of Theorems $4.9$ in \cite{ozawapopa}, B in \cite{ozawapopaII} and $3.3$ in \cite{houdayer7}.

\begin{theo}\label{step}
Let $G$ be an amenable group together with an orthogonal representation $\pi : G \to \mathcal{O}(H_\R)$. Let $M = \Gamma(H_\R, G, \pi)''$. Let $P \subset M$ be an amenable subalgebra such that $P \npreceq_M L(G)$. Then $\mathcal{N}_M(P)''$ is amenable.
\end{theo}

\begin{proof}
The proof is conceptually similar to the one of Theorem 4.9 in \cite{ozawapopa} under weaker assumptions: the malleable deformation $(\alpha_t)$ defined on $M = \Gamma(H_\R, G, \pi)''$ is not  assumed to be {\textquotedblleft compact over $L(G)$\textquotedblright} and the bimodule $L^2(\widetilde{M}) \ominus L^2(M)$ is merely weakly contained in the coarse bimodule $L^2(M) \otimes L^2(M)$. To overcome these technical difficulties, we will use ideas from the proof of Theorem B in \cite{ozawapopaII}. Note that the symbol {\textquotedblleft Lim\textquotedblright} will be used for a state on $\ell^\infty(\N)$, or more generally on $\ell^\infty(I)$ with $I$ directed, which extends the ordinary limit.

Let $G$ be an amenable group and let $\pi : G \to \mathcal{O}(H_\R)$ be an orthogonal representation. Let $M = \Gamma(H_\R, G, \pi)''$. Let $P \subset M$ be an amenable von Neumann subalgebra such that $P \npreceq_M L(G)$. Since $M$ has the c.m.a.p., $P$ is weakly compact inside $M$. Then there exists a net $(\eta_n)$ of vectors in $L^2(P \bar{\otimes} \bar{P})_+$ such that
\begin{enumerate}
\item $\lim_n \|\eta_n - (v \otimes \bar{v})\eta_n\|_2 = 0$, $\forall v \in \mathcal{U}(P)$; 
\item $\lim_n \|\eta_n - \Ad(u \otimes \bar{u})\eta_n\|_2 = 0$, $\forall u \in \mathcal{N}_M(P)$;
\item $\langle (a \otimes 1)\eta_n, \eta_n\rangle = \tau(a) = \langle \eta_n, (1 \otimes \bar{a}) \eta_n \rangle$, $\forall a \in M, \forall n$.
\end{enumerate}
We consider $\eta_n \in L^2(M \bar{\otimes} \bar{M})_+$, and note that $(J \otimes \bar{J}) \eta_n = \eta_n$, where $J$ denotes the canonical anti-unitary on $L^2(M)$. We shall simply denote $\mathcal{N}_M(P)$ by $\mathcal{G}$.

Let $z \in \mathcal{Z}(\mathcal{G}' \cap M)$ be a non-zero projection. Since $P \npreceq_M L(G)$ and $z \in P' \cap M$, it follows that $Pz \npreceq_M L(G)$. Theorem \ref{intertwining} then yields that the deformation $(\alpha_t)$ does not converge uniformly on $(Pz)_1$. Since any selfadjoint element $x \in (Pz)_1$ can be written
\begin{equation*}
x = \frac12 \|x\|_\infty (u + u^*)
\end{equation*}
where $u \in \mathcal{U}(Pz)$, it follows that $(\alpha_t)$ does not converge uniformly on $\mathcal{U}(Pz)$ either. Combining this with the inequality $(\ref{trans})$ in Proposition \ref{transversality}, we get that there exist $0 < c < 1$, a sequence of positive reals $(t_k)$ and a sequence of unitaries $(u_k)$ in $\mathcal{U}(P)$ such that $\lim_{k} t_k = 0$ and $\| \alpha_{t_k}(u_k z) - (E_M \circ \alpha_{t_k})(u_k z) \|_2 \geq c \|z\|_2$, $\forall k \in \N$. Since $\|\alpha_{t_k}(u_k z)\|_2 = \|z\|_2$, by Pythagora's theorem, we obtain
\begin{equation}\label{key}
\|(E_M \circ \alpha_{t_k})(u_k z)\|_2 \leq \sqrt{1 - c^2} \|z\|_2, \forall k \in \N.
\end{equation}
Set $\delta = \frac{1 - \sqrt{1 - c^2}}{6} \|z\|_2$. Choose and fix $k_0 \in \N$ such that 
\begin{equation}\label{delta}
\|\alpha_{t_k}(z) - z\|_2 \leq \delta, \forall k \geq k_0.
\end{equation}

Define for any $n$ and any $k \geq k_0$,
\begin{eqnarray*}
\eta_n^k & = & (\alpha_{t_k} \otimes 1)(\eta_n) \in L^2(\widetilde{M}) \otimes L^2(\bar{M}) \\
\xi_n^k & = & (e_M\alpha_{t_k} \otimes 1)(\eta_n) \in L^2(M) \otimes L^2(\bar{M}) \\
\zeta_n^k & = & (e_M^\perp\alpha_{t_k} \otimes 1)(\eta_n) \in (L^2(\widetilde{M}) \ominus L^2(M)) \otimes L^2(\bar{M}).
\end{eqnarray*}
We observe that
\begin{equation}\label{norm2}
\|(x \otimes 1) \eta_n^k\|_2^2 = \tau(E_M(\alpha_{t_k}^{-1}(x^*x))) = \|x\|_2^2, \forall x \in \widetilde{M}.
\end{equation}
As in the proof of Theorem $4.9$ in \cite{ozawapopa}, noticing that $L^2(\widetilde{M}) \otimes L^2(\bar{M})$ is an $M \bar{\otimes} \bar{M}$-module and since $\eta_n^k = \xi_n^k + \zeta_n^k$, Equation \eqref{norm2} gives that for any $u \in \mathcal{G}$, and for any $k \geq k_0$,
\begin{eqnarray}\label{crucial}
\mathop{\Lim}_n \|[u \otimes \bar{u}, \zeta^k_n]\|_2 
& \leq & \mathop{\Lim}_n \|[u \otimes \bar{u}, \eta_n^k]\|_2 \\ \nonumber
& \leq & \mathop{\Lim}_n \|(\alpha_{t_k} \otimes 1)([u \otimes \bar{u}, \eta_n])\|_2 + 2 \|u - \alpha_{t_k}(u) \|_2 \\
& = & 2 \|u - \alpha_{t_k}(u)\|_2. \nonumber
\end{eqnarray}
Moreover, for any $x \in M$,
\begin{eqnarray*}
\| (x \otimes 1) \zeta^k_n \|_2 & = & \|(x \otimes 1) (e_M^\perp \otimes 1) \eta_n^k\|_2 \\
& = & \|(e_M^\perp \otimes 1) (x \otimes 1) \eta_n^k\|_2 \\
& \leq & \| (x \otimes 1) \eta_n^k\|_2  = \|x\|_2.
\end{eqnarray*}

\begin{claim}\label{claim1}
For any $k \geq k_0$, 
\begin{equation}\label{crucial2}
\mathop{\Lim}_n \|(z \otimes 1) \zeta_n^k\|_2 \geq \delta.
\end{equation}
\end{claim}

\begin{proof}[Proof of Claim $\ref{claim1}$]
We prove the claim by contradiction. Exactly as in the proof of Theorem 4.9 in \cite{ozawapopa}, noticing that $e_M z = z e_M$ (since $z \in M$) and $z u_k = u_k z$ (since $z \in \mathcal{Z}(\mathcal{G}' \cap M)$), and using $(\ref{delta})$ we  have 
\begin{eqnarray*}
&& \mathop{\Lim}_n \|(z \otimes 1)\eta_n^k - (e_M \alpha_{t_k}(u_k) z \otimes \bar{u}_k)\xi_n^k\|_2  \\
& \leq & \mathop{\Lim}_n \|(z \otimes 1)\eta_n^k - (e_M \alpha_{t_k}(u_k) z \otimes \bar{u}_k)\eta_n^k\|_2 + \mathop{\Lim}_n \|(z \otimes 1) \zeta_n^k\|_2 \\
& \leq & \mathop{\Lim}_n \|(z \otimes 1)\eta_n^k - (e_M z \alpha_{t_k}(u_k) \otimes \bar{u}_k)\eta_n^k\|_2 + \|[\alpha_{t_k}(u_k), z]\|_2 + \delta \\
& \leq & \mathop{\Lim}_n\|(z \otimes 1)\zeta_n^k\|_2 + \mathop{\Lim}_n \|\eta_n^k - (\alpha_{t_k}(u_k) \otimes \bar{u_k})\eta_n^k\|_2 \\
& & + 2\|z - \alpha_{t_k}(z)\|_2 + \delta \\
& \leq & \mathop{\Lim}_n \|(\alpha_{t_k} \otimes 1)(\eta_n - (u_k \otimes \bar{u}_k)\eta_n)\|_2 + 4\delta = 4 \delta.
\end{eqnarray*}
Thus, we would get
\begin{eqnarray*}
\|(E_M \circ \alpha_{t_k})(u_k z)\|_2 & \geq & \|(E_M \circ \alpha_{t_k})(u_k)z\|_2 - \|z - \alpha_{t_k}(z)\|_2 \\
& \geq & \mathop{\Lim}_n \|((E_M \circ \alpha_{t_k})(u_k)z \otimes \bar{u}_k)\eta_n^k\|_2 - \delta \\
& \geq & \mathop{\Lim}_n \|(e_M \otimes 1) ((E_M \circ \alpha_{t_k})(u_k) z \otimes \bar{u}_k)\eta_n^k\|_2 - \delta \\
& = & \mathop{\Lim}_n \|(e_M \alpha_{t_k}(u_k) z \otimes \bar{u}_k) \xi_n^k\|_2 -  \delta \\
& \geq & \mathop{\Lim}_n \|(z \otimes 1)\eta_n^k\|_2 -  5\delta \\
& = & \|z\|_2 -  5\delta > \sqrt{1 - c^2} \|z\|_2,
\end{eqnarray*}
which is a contradiction according to $(\ref{key})$. 
\end{proof}

We now use the techniques of the proof of Theorem B in \cite{ozawapopaII}. Define a state $\varphi^{z, k}$ on $\mathbf{B}(\mathcal{H}) \cap \rho(M^{\op})'$, where $\rho(M^{\op})$ is the right $M$-action on $\mathcal{H}$, by
\begin{equation*}
\varphi^{z, k}(x) = \mathop{\Lim}_n \frac{1}{\|\zeta_n^{z, k}\|_2^2} \langle (x \otimes 1) \zeta_n^{z, k}, \zeta_n^{z, k}\rangle,
\end{equation*}
where $\zeta_n^{z, k} = (z \otimes 1) \zeta_n^k$. Note that 
\begin{equation*}
\varphi^{z, k}(x) = \varphi^{z, k}(z x) = \varphi^{z, k}(x z), \forall x \in \mathbf{B}(\mathcal{H}) \cap \rho(M^{\op})'.
\end{equation*}

\begin{claim}\label{claim2}
Let $a \in \mathcal{G}''$. Then one has 
\begin{equation*}
\mathop{\Lim}_k | \varphi^{z, k} (a x - x a) | = 0,
\end{equation*}
uniformly for $x \in \mathbf{B}(\mathcal{H}) \cap \rho(M^{\op})'$ with $\|x\|_\infty \leq 1$.
\end{claim}

\begin{proof}[Proof of Claim $\ref{claim2}$]
For $u \in \mathcal{G}$, since $z \in \mathcal{Z}(\mathcal{G}' \cap M)$, one has
\begin{eqnarray*}
\mathop{\Lim}_n \|\zeta_n^{z, k} - (u \otimes \bar{u}) \zeta_n^{z, k} (u \otimes \bar{u})^*\|_2 & \leq  & \mathop{\Lim}_n \|\zeta_n^{k} - (u \otimes \bar{u}) \zeta_n^{k} (u \otimes \bar{u})^*\|_2 \\
& \leq & 2 \|u - \alpha_{t_k}(u)\|_2.
\end{eqnarray*}
For every $x \in \mathbf{B}(\mathcal{H}) \cap \rho(M^{\op})'$, one has
\begin{equation*}
\varphi^{z, k}(u^* x u) = \mathop{\Lim}_n \frac{1}{\|\zeta_n^{z, k}\|_2^2} \langle(x \otimes 1)(u \otimes \bar{u}) \zeta_n^{z, k} (u \otimes \bar{u})^*, (u \otimes \bar{u}) \zeta_n^{z, k} (u \otimes \bar{u})^*\rangle,
\end{equation*}
so that with $(\ref{crucial}) - (\ref{crucial2})$,
\begin{equation*}
|\varphi^{z, k}(u^* x u) - \varphi^{z, k}(x)| \leq \frac{4}{\delta^2} \|x\|_\infty \|u - \alpha_{t_k}(u)\|_2.
\end{equation*}
This implies that 
\begin{equation*}
\mathop{\Lim}_k |\varphi^{z, k}(a x - x a)| = 0,
\end{equation*}
for each $a \in \mbox{span }\mathcal{G}$ and uniformly for $x \in \mathbf{B}(\mathcal{H}) \cap \rho(M^{\op})'$ with $\|x\|_\infty \leq 1$. However, for any $a \in M$,
\begin{eqnarray*}
|\varphi^{z, k}(x a)| & = & \mathop{\Lim}_n \frac{1}{\|\zeta_n^{z, k}\|_2^2} |\langle (x \otimes 1)(a \otimes 1) \zeta_n^{z, k}, \zeta_n^{z, k}\rangle| \\
& \leq & \frac{1}{\delta^2} \|x\|_\infty \|z a\|_2 \\
& \leq & \frac{1}{\delta^2} \|x\|_\infty \| a\|_2,
\end{eqnarray*}
and likewise for $|\varphi^{z, k}(a x)|$. An application of Kaplansky density theorem does the job.
\end{proof}

To prove at last that $\mathcal{G}''$ is amenable, we will use (as in Theorem B in \cite{ozawapopaII}) Connes' criterion for finite amenable von Neumann algebras (see Theorem $5.1$ in \cite{connes76} for the type ${\rm II_1}$ case and Lemma 2.2 in \cite{haagerup83} for the general case). For any non-zero projection $z \in \mathcal{Z}(\mathcal{G}' \cap M)$ and any finite subset $F \subset \mathcal{U}(\mathcal{G}'')$, we need to show
\begin{equation*}
\| \sum_{u \in F} uz \otimes \overline{uz} \|_{M \bar{\otimes} \bar{M}} = |F|.
\end{equation*}

Let $z \in \mathcal{Z}(\mathcal{G}' \cap M)$ be a non-zero projection and let $F \subset \mathcal{U}(\mathcal{G}'')$ be a finite subset. Since the $M, M$-bimodule $\mathcal{H}$ is weakly contained in the coarse bimodule $L^2(M) \otimes L^2(M)$, let $\Psi : \mathbf{B}(L^2(M)) \to \mathbf{B}(\mathcal{H}) \cap \rho(M^{\op})'$ be the u.c.p. map which extends the left $M$-action on $\mathcal{H}$ (see Lemma \ref{weakcontainment}). Note that $M$ is contained in the multiplicative domain of $\Psi$. Define $\psi^{z, k} = \varphi^{z, k} \circ \Psi$ a state on $\mathbf{B}(L^2(M))$. Let $u \in \mathcal{G}''$. By Claim \ref{claim2}, one has
\begin{eqnarray*}
\mathop{\Lim}_k |\psi^{z, k}((uz)^* x (uz) - x)| & = & \mathop{\Lim}_k |\varphi^{z, k}( \Psi((uz)^* x (uz)) - \Psi(x))| \\
& = & \mathop{\Lim}_k |\varphi^{z, k}( (uz)^* \Psi(x) (uz) - \Psi(x))| \\
& = & \mathop{\Lim}_k |\varphi^{z, k}( u^* \Psi(x) u - \Psi(x))| = 0,
\end{eqnarray*}
uniformly for $x \in \mathbf{B}(L^2(M))$ with $\|x\|_\infty \leq 1$. By a standard recipe of the theory together with the Hahn-Banach separation theorem, we can find a net $(\mu^{z, k})$ of positive norm-one elements in $S_1(L^2(M))$ (trace-class operators on $L^2(M)$) such that
\begin{equation*}
\lim_k \|\mu^{z, k} - \Ad(uz)\mu^{z,k}\|_1= 0, \forall u \in \mathcal{U}(\mathcal{G}'').
\end{equation*}
Since the above is satisfied in particular for $u = 1$ and since $F \subset \mathcal{\mathcal{G}''}$ is finite, replacing $\mu^{z, k}$ by $z \mu^{z, k} z/\|z \mu^{z, k} z\|_1$ we may assume that $\mu^{z, k} \in S_1(L^2(M))$ satisfies $\mu^{z, k} \geq 0$, $z\mu^{z, k}z = \mu^{z, k}$, $\|\mu^{z, k}\|_1 = 1$ and 
\begin{equation*}
\lim_k \|\mu^{z, k} - \Ad(uz)\mu^{z,k}\|_1= 0, \forall u \in F.
\end{equation*}
Define now $\nu^{z, k} = (\mu^{z, k})^{1/2} \in S_2(L^2(M))$ (Hilbert-Schmidt operators on $L^2(M)$). The net $(\nu^{z, k})$ satisfies $z\nu^{z, k}z = \nu^{z, k}$, $\|\nu^{z, k}\|_2 = 1$ and
\begin{equation*}
\lim_k \|\nu^{z, k} - \Ad(uz)\nu^{z,k}\|_2= 0, \forall u \in F.
\end{equation*}
by Powers-St\o rmer inequality. With the identification 
\begin{equation*}
S_2(L^2(M)) = L^2(M) \otimes L^2(\bar{M})
\end{equation*}
as $M, M$-bimodules it follows that the $\ast$-representations of $M$ and $\bar{M}$ given by the left and right $M$-actions induce the spatial tensor norm. Thus,
\begin{eqnarray*}
|F|  & = & \| \sum_{u \in F} \nu^{z, k}\|_2 \\
& \leq & \lim_k \| \sum_{u \in F} (uz)\nu^{z, k}(uz)^*\|_2 + \lim_k \| \sum_{u \in F} \nu^{z, k} - (uz) \nu^{z, k} (uz)^*\|_2 \\
& \leq & \|\sum_{u \in F} uz \otimes \overline{uz}\|_{M \bar{\otimes} \bar{M}}.
\end{eqnarray*}
Since the other inequality is trivial, the proof is complete.
\end{proof}

\subsection{Proof of Theorem A}

We refer to Section \ref{examples} for the necessary background on spectral measures of unitary representations. Let's begin with a few easy observations first. Assume that $(N, \tau)$ is a finite von Neumann algebra with no amenable direct summand, i.e. $Nz$ is not amenable, $\forall z \in \mathcal{Z}(N)$, $z \neq 0$. Then for any non-zero projection $q \in N$, $qNq$ is non-amenable. Moreover, if $N$ has no amenable direct summand and $N \subset N_1$ is a unital inclusion of finite von Neumann algebras, then $N_1$ has no amenable direct summand either.

\begin{lem}\label{diffuse}
Let $G$ be a countable group together with an action $G \curvearrowright (N, \tau)$ on a finite von Neumann algebra. Write $M = N \rtimes G$ for the crossed product. Let $B \subset N$ be a diffuse subalgebra. Then $B \npreceq_M L(G)$.
\end{lem}

\begin{proof}
We denote by $(v_g)$ the canonical unitaries which generate $L(G) \subset N \rtimes G = M$. Let $B \subset N$ be a diffuse subalgebra. Let $(u_n)$ be a sequence of unitaries in $B$ such that $u_n \to 0$ weakly, as $n \to \infty$. Let $I, J \subset G$ be finite subsets and 
\begin{eqnarray*}
x & = & \sum_{g \in I} x_g v_g \\
y & = & \sum_{h \in J} y_h v_h,
\end{eqnarray*}
where $x_g, y_h \in N$. Then we have
\begin{equation*}
E_{L(G)}(x^* u_n y) = \sum_{(g, h) \in I \times J} \tau(x_g^* u_n y_h) v_g^* v_h.
\end{equation*}
In particular,
\begin{equation*}
\|E_{L(G)}(x^* u_n y)\|_2 \leq \sum_{(g, h) \in I \times J} |\tau(x^*_g u_n y_h)|.
\end{equation*}
Since $u_n \to 0$ weakly, as $n \to \infty$, we get $\lim_n \|E_{L(G)}(x^* u_n y)\|_2 = 0$. Finally, using Kaplansky density theorem, we obtain
\begin{equation*}
\lim_n \|E_{L(G)}(x^* u_n y)\|_2 = 0, \forall x, y \in M.
\end{equation*}
By $(3)$ of Theorem \ref{intertwining1}, it follows that $B \npreceq_M L(G)$.
\end{proof}

\begin{theo}[Theorem A]\label{nocartan}
Let $\pi : \Z  \to \mathcal{O}(H_\R)$ be an orthogonal representation such that the spectral measure of $\pi$ has no atoms. Then $M = \Gamma(H_\R, \Z, \pi)''$ is a non-amenable ${\rm II_1}$ factor and for any maximal abelian subalgebra $A \subset M$, $\mathcal{N}_M(A)''$ is an amenable von Neumann algebra.
\end{theo}

\begin{proof}
Since the spectral measure of $\pi : \Z \to \mathcal{U}(H)$ has no atoms, it follows that $\pi$ has no eigenvectors. So the representation $\mathcal{F}(\pi) : \Z \to \mathcal{U}(\mathcal{F}(H))$ has no eigenvectors either. Thus, the corresponding free Bogoljubov action $\sigma^\pi : \Z \curvearrowright \Gamma(H_\R)''$ is necessarily outer (see Theorem~\ref{thm:bogoOuter}) and then $M = \Gamma(H_\R, \Z, \pi)''$ is a ${\rm II_1}$ factor. Moreover, $L(\Z)$ is clearly a MASA in $M$.

We prove the result by contradiction. Assume that $A \subset M= \Gamma(H_\R, \Z, \pi)''$ is a MASA such that $\mathcal{N}_M(A)''$ is not amenable. Write $1 - z \in \mathcal{Z}(\mathcal{N}_M(A)'')$ for the maximal projection such that $\mathcal{N}_M(A)''(1 - z)$ is amenable. Then $z \neq 0$ and $\mathcal{N}_M(A)''z$ has no amenable direct summand. Notice that $z \in A' \cap M = A$ and 
\begin{equation*}
\mathcal{N}_M(A)''z = \mathcal{N}_{zMz}(Az)'',
\end{equation*}
by Lemma 3.5 in \cite{popamal1}. Moreover $Az \subset zMz$ is a MASA.

Since the action $\sigma^\pi : \Z \curvearrowright \Gamma(H_\R)''$ is outer, it follows that $\Gamma(H_\R)' \cap M = \C$. Thanks to Theorem 3.3 in \cite{masapopa}, we can find a diffuse abelian subalgebra $B \subset \Gamma(H_\R)''$ which is a MASA in $M$. Since $M$ is a ${\rm II_1}$ factor and $B$ is diffuse, there exist a projection $p \in B$ and a unitary $u \in \mathcal{U}(M)$ such that $p = u z u^*$. Define $\tilde{A} = u Az u^*$. Then $\tilde{A} \subset pMp$ is a MASA and $\mathcal{N}_{pMp}(\tilde{A})''$ has no amenable direct summand. Let $C = \tilde{A} \oplus B(1 - p) \subset M$. Note that $C \subset M$ is still a MASA. Since $\mathcal{N}_M(C)''$ is not amenable and $C \subset M$ is weakly compact, Theorem \ref{step} yields $C \preceq_M L(\Z)$. Since $L(\Z)$ is a MASA, if we apply Theorem A.1 of \cite{popa2001}, we obtain $v \in M$ a nonzero partial isometry such that $v^*v \in C' \cap M = C$, $q = vv^* \in L(\Z)$ and $v C v^* \subset L(\Z)q$. Since $C \subset M$ is also a MASA, we get $v C v^* = L(\Z)q$. Note that $v p v^* \neq 0$, because otherwise we would have $v B v^* = L(\Z)q$ and this would imply that $B \preceq_M L(\Z)$, a contradiction according to Lemma \ref{diffuse}. Thus, with $q' = v p v^*$ we obtain $v \tilde{A} v^* = L(\Z)q'$. Consequently $\mathcal{N}_{q' M q'}(L(\Z) q')''$ is not amenable.

However, as $L(\Z), L(\Z)$-bimodules we have the following isomorphism
\begin{equation*}
L^2(M) \cong \bigoplus_{n \geq 0} K^{(n)}_{\pi},
\end{equation*}
where $K^{(n)}_{\pi} = H^{\otimes n} \otimes \ell^2(\Z)$ (see Section \ref{preliminaries}). Since the spectral measure of $\pi$ has no atoms, it follows that $L(\Z) \subset M$ is a singular MASA, i.e. $\mathcal{N}_M(L(\Z))'' = L(\Z)$, and {\em a fortiori} $\mathcal{N}_{q' M q'}(L(\Z) q')'' = L(\Z) q'$ (by Lemma 3.5 in \cite{popamal1}). We have reached a contradiction.
\end{proof}

\subsection{Proof of Theorem B} 

\begin{theo}[Theorem B]\label{stronglysolid}
Let $\pi : \Z  \to \mathcal{O}(H_\R)$ be a mixing orthogonal representation. Then the non-amenable ${\rm II_1}$ factor $M = \Gamma(H_\R, \Z, \pi)''$ is strongly solid.
\end{theo}

\begin{proof}
Since the representation $\pi : \Z \to \mathcal{O}(H_\R)$ is mixing, it has no eigenvectors. So the representation $\mathcal{F}(\pi) : \Z \to \mathcal{U}(\mathcal{F}(H))$ has no eigenvectors either. Thus, the free Bogoljubov action $\sigma^\pi : \Z \curvearrowright \Gamma(H_\R)''$ is necessarily outer (see Theorem~\ref{thm:bogoOuter}) and then $M = \Gamma(H_\R, \Z, \pi)''$ is a ${\rm II_1}$ factor.

Let $P \subset M$ be a diffuse amenable von Neumann subalgebra. By contradiction assume that $ \mathcal{N}_M(P)''$ is not amenable. Write $1 - z \in \mathcal{Z}(\mathcal{N}_M(P)'')$ for the maximal projection such that $\mathcal{N}_M(P)''(1 - z)$ is amenable. Then $z \neq 0$ and $\mathcal{N}_M(P)''z$ has no amenable direct summand. Notice that 
\begin{equation*}
\mathcal{N}_M(P)''z \subset \mathcal{N}_{zMz}(Pz)''.
\end{equation*}
Since this is a unital inclusion (with unit $z$), $\mathcal{N}_{zMz}(Pz)''$ has no amenable direct summand either. Let $A \subset \Gamma(H_\R)''$ be a diffuse abelian subalgebra. Since $M$ is a ${\rm II_1}$ factor and $A$ is diffuse, there exist a projection $q \in A$ and a unitary $u \in \mathcal{U}(M)$ such that $q = u z u^*$. Define $Q = u Pz u^*$. Then $Q \subset qMq$ is diffuse, amenable and $\mathcal{N}_{qMq}(Q)''$ has no amenable direct summand. Let $B = Q \oplus A(1 - q) \subset M$. Note that $B \subset M$ is a unital diffuse amenable subalgebra. Since $\mathcal{N}_M(B)''$ is not amenable and $B \subset M$ is weakly compact, Theorem \ref{step} yields $B \preceq_M L(\Z)$. 

Thus, there exists $n \geq 1$, a non-zero partial isometry $v \in \mathbf{M}_{1, n}(\C) \otimes M$ and a (possibly non-unital) $\ast$-homomorphism $\psi : B \to L(\Z)^n$ such that $x v = v \psi(x)$, $\forall x \in B$. Observe that $q v \neq 0$, because otherwise we would have $vv^* \leq 1 - q$ and $x v  = v \psi(x)$, $\forall x \in A(1 - q)$. This would mean that $A(1 - q) \preceq_M L(\Z)$ and so $A \preceq_M L(\Z)$, which is a contradiction according to Lemma \ref{diffuse}. Write $q v = w |q v|$ for the polar decomposition of $qv$. It follows that $w \in \mathbf{M}_{1, n}(\C) \otimes M$ is a non-zero partial isometry such that $x w = w \psi(x)$, $\forall x \in Q$. This means exactly that $Q \preceq_M L(\Z)$. Note that $ww^* \in Q' \cap qMq \subset \mathcal{N}_{qMq}(Q)''$ and $w^*w \in \psi(Q)' \cap \psi(q) M^n \psi(q)$.

Since the $\tau$-preserving action $\Z \curvearrowright \Gamma(H_\R)''$ is mixing by assumption and $\psi(Q) \subset \psi(q)L(\Z)^n\psi(q)$ is diffuse, it follows from Theorem $3.1$ in \cite{popamal1} (see also Theorem D.4 in \cite{vaesbern}) that $w^*w \in \psi(q)L(\Z)^n\psi(q)$, so that we may assume $w^*w = \psi(q)$. Note that $w^* Q w = \psi(Q)$. Moreover since $\psi(Q)$ is diffuse, Theorem 3.1 in \cite{popamal1} yields that the quasi-normalizer of $\psi(Q)$ inside $\psi(q)M^n\psi(q)$ is contained in $\psi(q) L(\Z)^n \psi(q)$. In particular, we get
\begin{equation*}
\Ad(w^*)(ww^*\mathcal{N}_{qMq}(Q)'' ww^*) \subset  \psi(q) L(\Z)^n \psi(q).
\end{equation*}
Note that $\Ad(w^*) : ww^* M ww^* \to w^*w M^n w^*w$ is a $\ast$-isomorphism. Since $\psi(q)L(\Z)^n \psi(q)$ is amenable and $ww^*\mathcal{N}_{qMq}(Q)'' ww^*$ is non-amenable, we finally get a contradiction, which finishes the proof.
\end{proof}

The above theorem is still true for any amenable group $G$ (instead of $\Z$), and any mixing orthogonal representation $\pi : G \to \mathcal{O}(H_\R)$ such that the corresponding Bogoljubov action $\sigma^\pi : G \curvearrowright \Gamma(H_\R)''$ is properly outer, i.e. $\sigma^\pi_g$ is outer, for any $g \neq e$.

\section{New examples of strongly solid ${\rm II_1}$ factors}\label{examples}

\subsection{Spectral measures and unitary representations} Let $H$ be a separable complex Hilbert space. Let $G$ be a locally compact second countable (l.c.s.c.) abelian group together with $\pi : G \to \mathcal{U}(H)$ a $\ast$-strongly continuous unitary representation. Denote by $\widehat{G}$ the dual of $G$. It follows that $C^*(G) \cong C_0(\widehat{G})$ and $\pi$ gives rise to a $\ast$-representation $\sigma : C_0(\widehat{G}) \to \mathbf{B}(H)$ such that $\sigma(f_g) = \pi(g)$, for every $g \in G$, where $f_g(\chi) = \chi(g)$, $\forall \chi \in \widehat{G}$.

Recall that for any unit vector $\xi \in H$, there exists a unique probability measure on $\mu_\xi$ on $\widehat{G}$ such that 
\begin{equation*}
\int_{\widehat{G}} f \, d\mu_\xi = \langle \sigma(f) \xi, \xi\rangle.
\end{equation*}
Note that the formula makes sense for every bounded Borel function $f$ on $\widehat{G}$.

\begin{df}
Let $G$ be a l.c.s.c. abelian group together with $\pi : G \to \mathcal{U}(H)$ a $\ast$-strongly continuous  unitary representation. The {\it spectral measure} $\mathcal{C}_\pi$ of the unitary representation $\pi$ is defined as the measure class on $\widehat{G}$ generated by all the probability measures $\mu_\xi$, for $\xi \in H$, $\|\xi\| = 1$.
\end{df}
Recall that the {\em support} of a measure is the (closed) subset of all points for which every neighborhood has positive measure. The spectral measure $\mathcal{C}_\pi$ is said to be {\it singular} if for all the probability measures $\mu$ in $\mathcal{C}_\pi$, the support of $\mu$ has $0$ Haar measure. From now on, we will only consider the cases when $G = \Z$ or $\R$.

We identify the Pontryagin dual of $\R$ with $\R$ by the pairing $\R \times \R \ni (x, y) \mapsto e^{2\pi i x y}$. Define
\begin{eqnarray*}
p: \R & \to & \mathbf{T} = \R/\Z \\
x & \mapsto & x + \Z
\end{eqnarray*}
the canonical projection. For $\mu$ a probability measure on $\R$, the push-forward measure of $\mu$ on $\mathbf{T}$ is defined by $(p_\ast\mu)(A) = \mu(p^{-1}(A))=\mu(A+\Z)$, $\forall A \subset \mathbf{T}$ Borel subset. The convolution product is denoted by $\ast$. We shall write
\begin{equation*}
\mu^{\ast k} = \mu \ast \cdots \ast \mu
\end{equation*}
for the $k$-fold convolution product.

\begin{lem}\label{pushforward}
Let $\mu$ be a probability measure on $\R$. Write $\nu = p_\ast \mu$.
\begin{enumerate}
\item If $\mu$ is singular, then $\nu$ is singular.
\item For any $k \geq 1$, $(p_\ast \mu)^{\ast k}$ and $p_\ast(\mu^{\ast k})$ are absolutely continuous to each other.
\end{enumerate}
\end{lem}

\begin{proof}
Denote by $\lambda$ the Lebesgue measure on $\R$. We may identify $(\mathbf{T}, \Haar)$ with $([0, 1], \lambda)$ as probability spaces. We use the notation $\mu_1 \sim \mu_2$ for two measures absolutely continuous to each other.

$(1)$ Assume that $\mu$ is singular. Write $K$ for the support of $\mu$ and $K_n = K \cap [n, n+ 1[$. Clearly, $\supp(\nu) \subset p(K)$. We have
\begin{eqnarray*}
\Haar(p(K)) & \leq & \sum_{n \in \Z} \Haar(p(K_n)) \\
& = & \sum_{n \in \Z} \lambda(K_n) = 0.
\end{eqnarray*}
Thus $\Haar(\supp(\nu)) = 0$ and $\nu$ is singular.

$(2)$ Under the previous identification, we have for any $A \subset \mathbf{T}$ Borel subset
\begin{eqnarray*}
\nu(B) & = & \mu(B + \Z) \\
& = &\sum_{n \in \Z} (\mu \ast \delta_n)(B).
\end{eqnarray*}
Thus for any $k \geq 1$, we have
\begin{eqnarray*}
\nu^{\ast k} & = & \left( \sum_{n \in \Z} \mu \ast \delta_n \right)^{\ast k} \\
& \sim & \sum_{n \in \Z} \left( \sum \mu^{\ast k} \ast \delta_n \right) \\
& \sim & \sum_{n \in \Z} \mu^{\ast k} \ast \delta_n.
\end{eqnarray*}
Consequently $(p_\ast \mu)^{\ast k} \sim p_\ast(\mu^{\ast k})$.
\end{proof}

\subsection{Examples of strongly solid ${\rm II_1}$ factors}

Erd\"os showed in \cite{erdos} that the symmetric probability measure $\mu_\theta$ on $\R$, with $\theta = 5/2$, obtained as the weak limit of
\begin{equation*}
\left( \frac12 \delta_{-\theta^{-1}} + \frac12 \delta_{\theta^{-1}} \right) \ast \cdots \ast \left( \frac12 \delta_{-\theta^{-n}} + \frac12 \delta_{\theta^{-n}} \right) 
\end{equation*}
is singular w.r.t. the Lebesgue measure $\lambda$ and has a Fourier Transform
\begin{equation*}
\widetilde{\mu}_\theta(t) = \prod_{n \geq 1} \cos\left(\frac{t}{\theta^n}\right)
\end{equation*}
which vanishes at infinity, i.e. $\widetilde{\mu}(t) \to 0$, as $|t| \to \infty$.

\begin{exam}\label{singularmeasure}
Modifying the measure $\mu_\theta$, Antoniou \& Shkarin (see Theorem $2.5, {\rm v}$ in \cite{antoniou}) constructed an example of a symmetric probability measure $\mu$ on $\R$ such that:
\begin{enumerate}
\item The Fourier Transform of $\mu$ vanishes at infinity, i.e. $\widetilde{\mu}(t) \to 0$, as $|t| \to \infty$.
\item For any $n \geq 1$, the $n$-fold convolution product $\mu^{\ast n}$ is singular w.r.t. the Lebesgue measure $\lambda$.
\end{enumerate}
\end{exam}
Let $\mu$ be a symmetric probability measure on $\R$ as in Example \ref{singularmeasure} and consider $\nu = p_\ast \mu$ the push-forward measure on the torus $\mathbf{T}$. Since $\mu(X) = \mu(-X)$, for any Borel set $X \subset \R$, it follows that $\nu(A) = \nu(\overline{A})$, for any Borel set $A \subset \mathbf{T}$, where $\overline{A} = \{\bar{z} : z \in A\}$.

Let $\pi^\nu : \Z \to \mathcal{U}(L^2(\mathbf{T}, \nu))$ be the unitary representation defined by $(\pi^\nu_n f)(z) = z^n f(z)$, $\forall f \in L^2(\mathbf{T}, \nu)$, $\forall n \in \Z$. Note that moreover
\begin{equation*}
H_\R^\nu = \left\{ f \in L^2(\mathbf{T}, \nu) : \overline{f(z)} = f(\bar{z}), \forall z \in \mathbf{T} \right\}
\end{equation*}
is a real subspace of $L^2(\mathbf{T}, \nu)$ invariant under $\pi^\nu$. Indeed, for all $f, g \in H^\nu_\R$, 
\begin{eqnarray*}
\langle f, g \rangle & = & \int_{\mathbf{T}} f(z) \overline{g(z)} \, d\nu(z) \\
& = & \int_{\mathbf{T}} \overline{f(\bar{z})} g(\bar{z}) \, d\nu(z) \\
& = & \int_{\mathbf{T}} \overline{f(\bar{z})} g(\bar{z}) \, d\nu(\bar{z}) \\
& = & \int_{\mathbf{T}} \overline{f(z)} g(z) \, d\nu(z) \\
& = & \overline{\langle f, g \rangle}.
\end{eqnarray*}
By assumption and using Lemma \ref{pushforward}, it follows that:
\begin{enumerate}
\item The unitary representation $\pi^\nu : \Z \to \mathcal{U}(L^2(\mathbf{T}, \nu))$ is mixing.
\item The spectral measure of $\bigoplus_{n \geq 1} (\pi^\nu)^{\otimes n}$ is singular.
\end{enumerate}

Consider now the non-amenable ${\rm II_1}$ factor $M = \Gamma(H_\R^\nu, \Z, \pi^\nu)''$. Let $A = L(\Z)$. Since $\pi^\nu$ is mixing, $A$ is maximal abelian in $M$ and {\em singular}, i.e. $\mathcal{N}_M(A)'' = A$. Since the spectral measure of the unitary representation $\bigoplus_{n \geq 1} (\pi^\nu)^{\otimes n}$ is singular and because of the $A, A$-bimodule isomorphism 
\begin{equation*}
L^2(M) \cong \bigoplus_{n \geq 0} K^{(n)}_{\pi^\nu},
\end{equation*}
where $K^{(n)}_{\pi^\nu} = L^2(\mathbf{T}, \nu)^{\otimes n} \otimes \ell^2(\Z)$ (see Section \ref{preliminaries}), it follows that the $A, A$-bimodule $L^2(M)$ is disjoint form the coarse bimodule $L^2(A) \otimes L^2(A)$. Combining Voiculescu's result (see Corollary 7.6 in \cite{voiculescu96}) and the second-named author's result (see Proposition $9.2$ in \cite{shlya99}), it follows that the non-amenable ${\rm II_1}$ factor $M$ is not isomorphic to any interpolated free group factor $L(\F_t)$, $1 < t \leq \infty$. Moreover, our Theorem \ref{stronglysolid} yields that $M$ is strongly solid, hence has no Cartan subalgebra.

\begin{theo}[Corollary B]\label{singular-ssolid}
The ${\rm II_1}$ factor $M = \Gamma(H_\R^\nu, \Z, \pi^\nu)''$ is strongly solid, hence has no Cartan subalgebra. Nevertheless, for the maximal abelian subalgebra $A = L(\Z)$, the $A, A$-bimodule $L^2(M)$ is disjoint from the coarse bimodule $L^2(A) \otimes L^2(A)$. Thus, $M$ is never isomorphic to an interpolated free group factor.
\end{theo}

\begin{rem}
For $\theta = 3$, $\mu_\theta$ is the Cantor-Lebesgue measure on the ternary Cantor set. If we set $\nu = p_\ast \mu_\theta$, we get that for any $n \geq 1$, the $n$-fold convolution product $\nu^{\ast n}$ is singular w.r.t. the Lebesgue measure $\lambda$. In that case, the ${\rm II_1}$ factor $M = \Gamma(H_\R^\nu, \Z, \pi^\nu)''$ has no Cartan subalgebras and is not isomorphic to any interpolated free group factor (Corollary A).
\end{rem}

\subsection{Bimodule decompositions over MASAs.}
Recall that if $\mu$ is a probability measure on $[0,1]\times [0,1]$ so that its push-forwards by the projection maps onto the two copies of $[0,1]$ are Lebesgue absolutely continuous, then $L^2([0,1]\times [0,1],\mu)$ can be regarded as an $L^\infty[0,1]$, $L^\infty[0,1]$-bimodule via the action
\begin{multline*}
(f_1 \cdot \xi \cdot f_2)(x,y) = f_1 (x) \xi(x,y) f_2 (y),
\\ x,y\in [0,1],\quad f_j\in L^\infty[0,1], \quad \xi\in L^\infty ([0,1]\times[0,1],\mu).\end{multline*}

For a von Neumann algebra $M$, consider the collection $\mathcal{C}(M)$ of measure classes $[\mu]$ on $[0,1]\times [0,1]$ with the property that there exists a MASA $L^\infty [0,1]\cong A\subset M$ so that $L^2(M)$, when regarded as an $A,A$-bimodule, contains a copy of $L^2([0,1]^2,\mu)$.  
Also let $\mathcal{D}(M)$ be the collection of all measure classes $[\mu]$ so that for {\em every} MASA $L^\infty[0,1]\cong A\subset M$, $L^2(M)$ contains a sub-bimodule of $L^2([0,1]^2,\mu)$.  Clearly, $\mathcal{C}\supset \mathcal{D}$.  

Then (as is well known) $M$ has a Cartan subalgebra if and only if  $\mathcal{C}(M)$ contains an $r$-discrete measure class (i.e., a measure class $[\mu]$ for which $\mu(B) = \int \mu_t (B) dt$ and $\mu_t$ are a.e. discrete).

Voiculescu in \cite{voiculescu96} proved that $\mathcal{D}(L(\F_n)) \ni \{\textrm{Lebesgue Measure}\}$.

It thus remained open whether every II$_1$ factor $N$ must either contain a Cartan subalgebra, or satisfy that  $\mathcal{D}(N) \ni \{\textrm{Lebesgue Measure}\}$.  Our main example $M= \Gamma(H_\R^\nu, \Z, \pi^\nu)''$ answers this question in the negative, as $\mathcal{D}(M)$ does not contain Lebesgue measure and yet $M$ has no Cartan subalgebra.

\section{Outerness of free Bogoljubov actions}

Although we do not need the following result in the rest of the paper, we record the following observation, which is well-known to the experts and is most likely folklore (although we could not find a precise reference).

\begin{theo}\label{thm:bogoOuter} Let $G$ be a countable group, and let  $\pi: G \to \mathcal{O}(H_\R)$ be a $\ast$-strongly continuous orthogonal
representation of $G$ on a real Hilbert space $H_\R$.  Then $\sigma_g^\pi$ is inner iff $\pi_g=1$. 
In particular, if $\pi_g \neq 1$ for any $g\neq e$, the Bogoljubov action $\sigma^\pi$ of $G$ on $\Gamma(H_\R)''$  is outer.
\end{theo}
\begin{proof}
Let $g$ be an element of $G$ so that $\pi_g \neq 1$, 
and let $\alpha = \sigma_g^\pi$ acting on $M=\Gamma(H_\R)''$. Let $T=\pi_g$. We may assume without loss of generality that $H_\R$ has dimension at least $2$, so that $M$ is a factor (otherwise, $M$ is abelian, and any non-trivial $T$ gives rise to an outer transformation).

Suppose for a contradiction that $\alpha = \Ad(u)$ for some unitary $u\in M$.  Then for any $x\in M$,
$$
\alpha(x) = u x u^*
$$ and so $\alpha (u) = u$.

Let $H = H_\R \otimes_\R \C$ be the complexification of $H_\R$.  We continue to denote the complexification of $T$ by the same letter.
Let $H^a \subset H$ be the closed linear span of eigenvectors of $T$, $H^a_\R  = H^a \cap H_\R$ be its real part. Then $N=\Gamma(H^a_\R)''\subset \Gamma(H_\R)''=M$.  Moreover, it is clear from the Fock space decomposition of $L^2(M)$ that any eigenvectors for $\alpha$ must lie in $L^2(N)$, so $u\in N$.  Thus we may, without loss of generality, assume that $N=M$ and that eigenvectors of $T$ densely span $H$.

Thus we may assume that 
\begin{equation*}
H_\R = \R^n \oplus \bigoplus_{k\in J} H^k_\R,
\end{equation*}
where $n\in \{0,1,\dots,+\infty\}$,
each $H^k_\R\cong \R^2$ and $T$ acts trivially on $\R^n$ and acts on $H^k_\R$ by a rotation of period $2\pi /\log \lambda_k$.  If we denote by $h_k, g_k$ an orthonormal basis for $H^k_\R$ and we set $c_k = s(h_k) + i s(g_k)\in M$, then $M\cong L(\F_n) * W^*(c_k : k\in J)$, and $\alpha = \id * \beta$ where
$\beta (c_j ) = \exp(2\pi i \lambda_j) c_j$.   Let $c_j = u_j b_j$ be the polar decomposition of $c_j$; thus $\beta(u_j)=\exp(2\pi i \lambda_j) u_j$ and $\beta(b_j)=b_j$.
By \cite{DVV:circular}, $b_j$ and $u_j$ are freely independent and $W^*(b_k:k\in J) \cong W^*(u_k : k\in J) \cong L(\F_{2|J|})$.  It follows that $M\cong L(\F_n) * W^*(b_k : k\in J) * W^*(u_k : k\in J) \cong L(\F_{n+|J|}) * L(\F_{|J|}) = N * P$ in such a way that $\alpha$ corresponds to the action $\id * \gamma$ where $\gamma : P\to P = W^*(u_k : k\in J)$ is given by $\gamma(u_k)=\exp(2\pi i \lambda_k) u_k$.

Since by assumption $T$ is non-trivial, $|J|\geq 1$ and also $|J| + n \geq 1$.  Thus if $\alpha (x) = uxu^*$ for all $x\in M$, then $u$ must commute with  $N\subset  N * P\cong M$.  But $N'\cap M = N'\cap N = \mathcal{Z}(N)$ (e.g. because as an $N$,$N$-bimodule, $L^2(M) = L^2(N) \oplus (\textrm{a multiple of coarse $N$,$N$-bimodule})$), so $u\in \mathcal{Z}(N)$.  But then $u P u^* =\alpha(P)\subset P$, which is easily seen to be impossible by using the free product decomposition of $L^2(M)$ in terms of $L^2(N)$ and $L^2(P)$, unless $u=\tau(u)$.  But this is impossible, since $\alpha(s(h)) = s(Th)$ is a non-trivial automorphism.
\end{proof}

\section{Free Krieger algebras}

Let $\nu$ be a probability measure on the torus $\mathbf{T}$. Note that $\nu$ gives rise to unital completely positive map $\eta : A \to A$, ($A = L^\infty(\mathbf{T})$), determined by 
$$\eta(f)(x) = \int f(x-y) d\nu(y) = (f * \nu)(x), \forall f \in C(\mathbf{T}).$$
It is not hard to see that the von Neumann algebra $M=\Gamma(H_\R^\nu, \mathbf{Z},\pi^\nu)'' \cong \Phi (A,\eta)$ in the notation of \cite{shlya99}, i.e., it is an example of a von Neumann algebra generated by an $A$-valued semicircular system with covariance $\eta$ (these were called ``free Krieger algebras'' in \cite{shlya99}, following the analogy between the operation $A\mapsto \Phi(A,\eta)$ and the crossed product operation $A\mapsto A\rtimes_\sigma \mathbf{Z}$). 

 As we have seen, $M$ has both the c.m.a.p. and the Haagerup property, and thus for this specific choice of $\eta$, $\Phi(A,\eta)$ has these properties.  
 
We point out that in general (even for abelian $A$), $\Phi(A,\eta)$ may fail to have the Haagerup property for other choices of the completely positive maps $\eta$.  It is an interesting question to determine exactly when $\Phi(A,\eta)$ has this property (and/or c.m.a.p.) as a condition on the completely-positive map $\eta : A\to A$, $A\cong L^\infty[0,1]$.  It is likely that the techniques of the present paper would then apply to give solidity of $\Phi(A,\eta)$.

\begin{prop}  
There exists a choice of $\eta : A\to A$, $A\cong L^\infty [0,1]$, so that $\Phi(A,\eta)$ does not have the Haagerup property and is not weakly amenable, i.e. $\Lambda_{\cb}(\Phi(A,\eta)) = \infty$.
\end{prop}

\begin{proof}
Let $\alpha$ be an action of a free group $\F_2$ on $A\cong L^\infty[0,1]$ so that $M=A\rtimes_\alpha \F_2$ does not have the Haagerup property and is not weakly amenable (one could take, for example, an action measure equivalent to the action of $\SL(2, \mathbf{Z})$ on $A=L(\mathbf{Z}^2)$; the crossed product in this case has relative property (T) and does not have the Haagerup property  \cite{popa2001}.  Moreover it is not weakly amenable, i.e. $\Lambda_{\cb}(M) = \infty$ (see \cite{dorofaeff}). Denote the two automorphisms of $A$ corresponding to the actions of the two generators of $\F_2$ by $\alpha_1$, $\alpha_2$, and let $\eta_j = \alpha_j + \alpha_j^{-1}$, $\eta = \eta_1+\eta_2$.  

Let $\sigma$ be the free shift action of $\mathbf{Z}$ on $\F_\infty$.  Then by \cite{shlya99}, 
\begin{multline}
\Phi(A,\eta) \cong \Phi(A,\eta_1) *_A \Phi(A,\eta_2) \cong \\
\left((A \bar{\otimes} L(\F_\infty)) \rtimes_{\alpha_1 \otimes \sigma} \mathbf{Z}\right)
*_A
\left((A \bar{\otimes} L(\F_\infty)) \rtimes_{\alpha_2\otimes \sigma} \mathbf{Z}\right) \cong \\
(A \bar{\otimes} [L(\F_\infty)*L(\F_\infty)])\rtimes_{\alpha \otimes \sigma*\sigma}\F_2.
\end{multline}
Thus $\Phi(A,\eta)$ contains $M$ as a subalgebra.  Since the Haagerup property and the weak amenability are inherited by subalgebras, it follows that $\Phi(A,\eta)$ cannot have the Haagerup property and is not weakly amenable.
\end{proof}

%\addcontentsline{toc}{section}{Bibliography}
\bibliographystyle{plain}

\end{document}